\documentclass{birkjour}
 \theoremstyle{definition}
 
 \theoremstyle{remark}

 \numberwithin{equation}{section}

\usepackage{amsmath}
\usepackage{amssymb}
\usepackage{amsthm}
\usepackage{enumitem}

\setlist[enumerate,1]{label = \alph*),ref = \alph*}

\usepackage[colorinlistoftodos,textsize=small,backgroundcolor=white,linecolor=black]{todonotes}

\newcounter{todoMCcounter}

\newtheorem{definition}{Definition}[section]
\newtheorem{example}[definition]{Example}

\newtheorem{theorem}[definition]{Theorem}
\newtheorem{lemma}[definition]{Lemma}

\newcommand{\R}{\mathbb{R}}                             
\newcommand{\Z}{\mathbb{Z}}                             

\DeclareMathOperator{\JSymbol}{\mathcal{J}}	            

\newcommand{\embed}{\hookrightarrow}                    

\DeclareMathOperator{\LspaceSymbol}{\mathbf{L}}

\newcommand{\LPhispace}[1][\DefaultPhifunction]{\LspaceSymbol^{{#1}}}

\DeclareMathOperator{\WspaceSymbol}{\mathbf{W}}

\newcommand{\WLPhispace}[1][\DefaultPhifunction] { {\WspaceSymbol^1}\LspaceSymbol^{{#1}} }

\newcommand{\norm}[1]{\|#1\|}                           

\newcommand{\LPhinorm}[3][\DefaultPhifunction]{\norm{#2}_{\LPhispace[{#1}]({#3})}}														                                          	

\newcommand{\WLPhinorm}[3][\DefaultPhifunction]{\norm{#2}_{\WLPhispace[{#1}]({#3})}}

\begin{document}

%
%
%
%
%
%
%
%
%

\title[Elliptic problem in anisotropic Orlicz-Sobolev space]
 {Quasilinear elliptic problem in \\anisotropic Orlicz-Sobolev space on 
\\unbounded domain}

\author{Karol Wro\'{n}ski}
\address{%
Department of Technical Physics and Applied Mathematics\\
Gda\'{n}sk University of Technology\\
Narutowicza 11/12, 80-233 Gda\'{n}sk, Poland}
\email{karwrons@pg.edu.pl}
 


\date{January 1, 2004}

\begin{abstract}
We study a quasilinear elliptic problem $-\text{div} (\nabla \Phi(\nabla 
u))+V(x)N'(u)=f(u)$ with anisotropic convex function $\Phi$ on whole $\R^n$. To 
prove existence of a nontrivial weak solution we use mountain pass theorem for 
a functional defined on anisotropic Orlicz-Sobolev space $\WLPhispace(\R^n)$. 
As the domain is unbounded we need to use Lions type lemma formulated for Young 
functions. Our assumptions broaden the class of considered 
functions $\Phi$ so our result generalizes earlier analogous results proved in 
isotropic setting.
\end{abstract}

\maketitle
\section{Introduction}

We consider the following anisotropic quasilinear problem 

\begin{equation} \tag{AQP} \label{eq:AQP}
-\text{div} (\nabla \Phi(\nabla u))+V(x)N'(u)=f(u),\quad \text{ where }
u\in\WLPhispace(\R^n),
\end{equation} 
$V:\R^n\to\R$ and $f, N:\R\to\R$ are continuous  functions. We assume that
$\Phi:\R^n\to\R$ is a differentiable, strictly convex, n-dimensional N function 
(G-function), $N$ and $f$ satisfy some growth conditions and 
$\WLPhispace(\R^n)$ is the (possibly anisotropic) Orlicz-Sobolev space.

For the reader's convenience, all assumptions are provided below. The necessary
definitions of the terms and symbols used here are provided in the next
chapter. We assume that:

\begin{equation}\label{eq:Phi_delta}\tag{$\Delta$}
	\text{ $\Phi$ satisfies conditions $\Delta_2$,$\nabla_2$ globally,}
\end{equation}

\begin{equation}\label{eq:int<infty}\tag{$\Phi_0$}
	\int_0\left(\frac{t}{\Phi_{\circ}(t)}\right)^{\frac{1}{n-1}}\,dt<\infty
\end{equation}

\begin{equation}\label{eq:Phi_infty_diverg}\tag{$\Phi_1$}
\int^{\infty}\left(\frac{t}{\Phi_{\circ}(t)}\right)^{\frac{1}{n-1}}\,dt=\infty
\end{equation}

\begin{equation}
\label{as:Phi<< Phinzdachem}\tag{$\Phi_2$}
\Phi\prec\prec 	\widehat{\Phi_n}
\end{equation}

\begin{equation}\label{N:growth}\tag{$N_1$}
	\text{$N$ is a differentiable N-function, such that }
	N\approx \Phi_{\circ}
\end{equation}

\begin{equation}\label{f:in0}\tag{$f_1$}
	\lim_{t\to 0}\frac{f(t)t}{N(t)}=0,
\end{equation}

\begin{equation}\label{f:ininfty}\tag{$f_2$}
	\lim_{t\to 
		\infty}\frac{f(t)t}{\Phi_n(t)}=0,
\end{equation}

\begin{equation}\label{f:takitrocheAR}\tag{$f_3$}
	0<\theta F(t)=\leq tf(t) \quad\text{for some 
		$\theta>s_{\Phi},\, \theta>s_N$ and all } t\in\R\backslash\{0\},
\end{equation}
where $F(t)=\int_0^tf(s)\,ds$ and the formula for the index $s_{\Phi}$ can be 
found in Definition 
\ref{def:indices},

\begin{equation}\label{as:V>0}\tag{$V_1$}
	0<V_0=\inf_{x\in\R^n}V(x),\text{ and $V$ is $1$-periodic},
\end{equation}

\begin{equation}\label{as:V-ZNperiodic}\tag{$V_2$}
	V(x+k)=V(x)\quad  \text{for all }x\in\R^n, k\in\Z^n.
\end{equation}
Assumptions on $\Phi$ and $N$ will be given in the next chapter together with 
the definition of the space $\WLPhispace(\R^n)$.

\begin{theorem}[Main Theorem] \label{thm:main} Suppose that $V$ and $f$ satisfy  
\eqref{f:in0}-\eqref{f:takitrocheAR},
	\eqref{as:V>0}-\eqref{as:V-ZNperiodic} and
 $\Phi$, $N$ satisfy assumptions \eqref{eq:int<infty} - \eqref{as:Phi<< 
Phinzdachem}, \eqref{eq:Phi_delta}, \eqref{N:growth}. Then the problem
\eqref{eq:AQP} has a nontrivial weak solution. \end{theorem}

The theorem will be proved using the mountain pass theorem applied for the 
functional 
$$\JSymbol(u)=\int_{\R^N}\Phi(\nabla u)+V(x)N(u)-F(u)\,dx.$$ 
A theorem of this 
type was proved in \cite{AlvFigGioSan14} where authors work 
with standard isotropic (i.e. depending only on the norm of a vector) Young 
function instead of function $\Phi$ depending on 
$n$ variables. Theorem \ref{thm:main} generalizes the analogous result from 
\cite{AlvFigGioSan14} because our assumptions broaden the class of considered 
functions $\Phi$. Our proof techniques are similar but it occurred
that for technical reasons the proof cannot be simply rewritten in  
anisotropic case. For similar problems with laplacian and p-laplacian see 
references in  
\cite{AlvFigGioSan14}.

As the problem is considered on unbounded domain we cannot use compact 
embeddings of Orlicz-Sobolev space into Orlicz space. Thus there is a 
problem with convergence of the minimizing sequence. The main tool  which  is 
used to deal with lack of compactness on unbounded domain is the Lions lemma 
originally formulated in \cite{Lio84}.
The formulation of this lemma and its  application for equation with Laplace  
operator can be found in \cite{Cos07}. 
Theorem similar to theorem \ref{thm:main}, but formulated for fractional 
Orlicz-Sobolev space, can be found in \cite{SilCarAlbBah21}. We use Lions type 
lemma in a version given in \cite{AlvFigGioSan14}.

To the best author's knowledge problem \eqref{eq:AQP} with anisotropic function  
$\Phi$ for an unbounded domain has not been considered before.  Examples of 
applications of isotropic  problem can be found for example in 
\cite{AlvFigGioSan14}. Anisotropy significantly widens these applications. More 
information about anisotropic operators in divergence form can be found in 
e.g. \cite{BarCia17} and  \cite{AlChlCiaZat19}. 

Techniques specific to our case differ from the cited ones mainly in the 
steps when we use lemmas \ref{lem:Phinzdachem} and \ref{lem:PhimodularNormcon}. 
Those lemmas 
turned out to be useful in most of the proof steps of the main theorem 
including the properties of the functional and many others. We also use 
specific assumption \eqref{as:Phi<< Phinzdachem} which in the isotropic case 
would be trivially satisfied for every Young function but in our anisotrpic 
setting is not trivial - see example \ref{ex:growth}. This assumption is crucial 
in lemma \ref{lem:Phinzdachem} and its consequences.

\section{Preliminaries}

Below we provide the necessary set of definitions and facts from the theory of 
anisotropic convex functions and Orlicz-Sobolew spaces. This part does not
contain proofs of some well known facts, the reader can find more information
about the subject in for example \cite{BarCia17}, \cite{Cia00} and especially in
\cite{Shap05} where a very general setting is presented.

A function $\Phi:\R^n\to[0,\infty)$  is 
called an n-dimensional 
Young function  if it is convex, even, $\Phi(0)=0$ and 
$\{v\in\R^n\colon\Phi(v)\leq t\}$ is a compact set containing
$0$ in its interior for every $t > 0$. The function $\Phi:\R^n\to[0,\infty)$ is 
called an n-dimensional N-function if, 
 in addition, $\Phi$ is finite–valued, vanishes only at 0, and
\[
\lim_{v\to0}\frac{\Phi(v)}{|v|}=0\quad \text{and}\quad  
\lim_{|v|\to\infty}\frac{\Phi(v)}{|v|}=\infty.
\]

The assumption \eqref{eq:Phi_delta} that $\Phi$ satisfies conditions
$\Delta_2$,$\nabla_2$ globally, means that there exist $K_1>2$, $K_2\geq 2$,
such that
$K_1\Phi(v)\leq\Phi(2v)\leq
K_2\Phi(v)$ for all $v\in \R^n$.

Below we will give definitions and list some facts about Orlicz spaces.
Such spaces were considered in many papers and books. Our setting is not the
most popular in literature because the domain may be whole $\R^n$ and we
use n-dimensional N-functions. However this situation also has been
considered, e.g. \cite{Shap05}.

An anisotropic Orlicz space $\LPhispace(\R^n)$ is  the set of all measurable
functions $v\colon\R^n\to\R^n$, such that the Luxemburg norm
\[
\LPhinorm{v}{\R^n}=\inf\left\{k>0: \int_{\R^n} \Phi\left(\frac{v}{k}\right)\, dt\leq 1
\right\}
\]	
 is finite and isotropic Orlicz space $\LPhispace[N](\R^n)$ as  the
 set of all measurable functions $u\colon\R^n\to\R$, such that the Luxemburg norm
 \[
 \LPhinorm[N]{u}{\R^n}=\inf\left\{k>0: \int_{\R^n} N\left(\frac{|u|}{k}\right)\, dt\leq 1
 \right\}
 \]	
 is finite. Note that we define Orlicz space for both scalar and vector 
valued functions.

Now we can define the Orlicz-Sobolev space as a set of weakly differentiable
functions on $\R^n$, such that  the norm 
$$\WLPhinorm{u}{\R^n}=\LPhinorm{\nabla u}{\R^n}+\LPhinorm[N]{u}{\R^n}$$ is finite.
Since $\Phi$ satisfies $\Delta_2$ and $\nabla_2$ conditions we have that 
$\WLPhispace$ is a reflexive Banach space. Another important property of 
functions satisfying $\Delta_2$ and $\nabla_2$ is that convergence in norm 
$\|\cdot\|_{\LPhispace}$ is equivalent to the convergence of modular, i.e. 
$\int\Phi(x_n-x)dx\to 0\Leftrightarrow \|x_n-x\|_{\LPhispace}\to 0$.

Denote by  $\Phi_{\circ}:[0,\infty)\to[0,\infty)$ the function obeying
\begin{equation}\label{eq:Phi_circ}
\lambda(\{v\in\R^n\colon\Phi_{\circ}(|v|)\leq 
t\})=\lambda(\{v\in\R^n\colon\Phi(v)\leq t\}) 
\quad \text{ for  } t\geq0,
\end{equation}
where $\lambda$ stands for Lebesgue measure. Function $\Phi_{\circ}$ is a
convex function as a consequence of the Brunn - Minkowski inequality. This fact
is stated in many papers but we could not find any straightforward proof for it
so below we write it as a lemma. We will also prove that $\Phi_{\circ}$
satisfies $\Delta_2$ and $\nabla_2$ globally.

\begin{lemma}
 $\Phi_{\circ}$ is convex and satisfies $\Delta_2$ and $\nabla_2$ globally.
\end{lemma}
\begin{proof}
We begin with an easy observation that by definition
$$\Phi_{\circ}(s)=\sup\{t: \lambda(\{v\in\R^n\colon\Phi(v)\leq
t\})\leq\omega_n s^n\}$$
where $\omega_n$ denotes the measure of unit ball. We define a set
$B_t=\{v\in\R^n\colon\Phi(v)\leq t\}$ and a function $\beta:
[0,\infty)\to[0,\infty)$, $\beta(t)=(\lambda(B_t))^{\frac1n}$. Notice that
$\beta(0)=0$, $\beta$ is increasing and continuous,
$\lim\limits_{t\to\infty}\beta(t)=\infty$. This properties of $\beta$ together
with the formula for $\Phi_{\circ}$ given above give us
$\Phi_{\circ}(s)=\beta^{-1}(\omega^{\frac1n}s)$. Now we will prove that $\beta$
is concave.

Chose $r,t\geq 0$,  $\varepsilon\in[0,1]$ and $s=(1-\varepsilon)r+\varepsilon\,
t$. By convexity of $\Phi$ we have
$$(1-\varepsilon)B_r+\varepsilon B_t=\{(1-\varepsilon)x+\varepsilon\, y\, :\,
\Phi(x)\leq r, \Phi(y)\leq t \}\subset B_s$$
because $\Phi((1-\varepsilon)x+\varepsilon\, y)\leq
(1-\varepsilon)\Phi(x)+\varepsilon\, \Phi(y)\leq s$ if $\Phi(x)\leq r$ and
$\Phi(y)\leq t$. By the Brunn-Minkowski inequality
\begin{align*}
\beta(s)&=(\lambda(B_s))^{\frac1n}\geq (\lambda((1-\varepsilon)B_r+\varepsilon
B_t))^{\frac1n}\\
&\geq (1-\varepsilon)(\lambda(B_r))^{\frac1n}+
\varepsilon(\lambda(B_t))^{\frac1n}=
(1-\varepsilon)\beta(r)+\varepsilon\, \beta(t).
\end{align*}
Thus $\beta$ is concave and $\Phi_{\circ}(s)=\beta^{-1}(\omega^{\frac1n}s)$ is
convex.

Now we will prove the $\Delta_2$ property of $\Phi_{\circ}$. As $\Phi$
satisfies $\Delta_2$ we have
\begin{align*}
 (\lambda(B_{K_2t}))^{\frac1n}&= (\lambda(\{v\in\R^n\colon
\frac1{K_2}\Phi(v)\leq t\})^{\frac1n}\geq (\lambda(\{v\in\R^n\colon \Phi(\frac
v2)\leq t\})^{\frac1n}\\
&= (\lambda(\{2v\in\R^n\colon \Phi(v)\leq t\})^{\frac1n}=
2\lambda(B_t)^{\frac1n}
\end{align*}
so $2\beta(t)\leq\beta(K_2t)$ and $\beta^{-1}(2\beta(t))\leq
K_2t$. For $y=\beta(t)$ we get $\beta^{-1}(2y)\leq K_2\beta^{-1}(y)$. The proof
of $\nabla_2$ property is analogous.
\end{proof}
	
As we assume \eqref{eq:int<infty} in the same way as in \cite{Cia00}, for
$s\geq 0 $ we can define $H:[0,\infty)
\to[0,\infty)$ by
\begin{equation}
	\label{eq:H}
	H(s)=\left(\int_0^s\left(\frac{t}{\Phi_{\circ}(t)}\right)^{\frac{1}{n-1}} 
\,dt\right)^{\frac{n-1}{n}}
\end{equation}
Define a Young function $\Phi_n:[0,\infty)\to[0,\infty]$ by 
\begin{equation}
\label{eq:phin}
	\Phi_n=\Phi_{\circ}\circ H^{-1},
\end{equation}
where $H^{-1}$ is the left-continuous inverse of $H$. 

\begin{example}[p. 51, \cite{BarCia17}]
	If $\Phi:\R^n\to[0,\infty)$, $\Phi(v)=\sum|v_i|^{p_i}$, then $\Phi_n(|v|) 
\approx|v|^{{\overline p}^*},$  $\overline{p}$ is a  
harmonic average of the powers i.e.  $\overline{p}=\frac{n}{\sum\frac{1}{p_i}}$ 
and 
	${{\overline p}^*}=\frac{n\overline{p}}{n-\overline{p}}.$
\end{example}

\begin{theorem}[thm. 1, \cite{Cia00}] \label{thm:embcian}There exists constant  
$K>0$ depending only on $n$, such that
	\begin{equation} \label{eq:poinc:int}
		\int_{\R^n}\Phi_n\left(\frac{|u(x)|}{K\left(\int_{R^n}\Phi(\nabla u)\,dy 
\right)^{1/n}}\right)\,dx\leq\int_{\R^n}\Phi(\nabla u)\,dx
	\end{equation}
	and
	\begin{equation} \label{eq:poinc:norm}
		\LPhinorm[\Phi_n]{u}{\R^n}\leq K\LPhinorm{\nabla u}{\R^n}
	\end{equation}
	for n-dimensional Young function $\Phi$ satisfying \eqref{eq:int<infty} and  
for every real-valued weakly differentiable function $u$ satisfying 
$|\{x\in\R^n\colon |u(x)|>t\}|<\infty$ for every $t>0$. 
\end{theorem}

Above theorem establishes embedding of $\WLPhispace$ into  Orlicz space  
$\LPhispace[\Phi_n]$. It is worth to notice that there are two 
important cases 
of this embedding. If we assume that \eqref{eq:Phi_infty_diverg} does not hold,
which means
\begin{equation}\label{eq:Phi_infty_converg}
\int^{\infty}\left(\frac{t}{\Phi_{\circ}(t)}\right)^{\frac{1}{n-1}}\,dt<\infty
\end{equation}
then $\Phi_n(t)=\infty$ for large $t$ and in fact $\LPhispace[\Phi_n]$ turns out 
to be $\mathbf{L^{\infty}}$. If on the other hand we assume
\eqref{eq:Phi_infty_diverg}
then $\Phi_n$ is a scalar valued function and the space
$\LPhispace[\Phi_n]$ is larger than $\mathbf{L^{\infty}}$. This case is more 
complicated and 
we will restrict ourselves only to this case.

\begin{definition}\label{def:growth_relation}
We say that G-functions are in relation $A\prec B$ if for some constant $C$ and
large $v$ we have $A(v)\leq B(C\, v)$. Relation $A\approx
B$ means that both $A\prec B$ and $B\prec A$. Relation $A\prec\prec B$ is
satisfied when
\[\lim\limits_{x\to\infty}\frac{B(cx)}{A(x)}=\infty\quad\text{for all }c>0.\]
\end{definition}

Recall that $A\prec B$ implies $\LPhispace[B]\subset\LPhispace[A]$, and if
$A\prec\prec B$ the embedding is also compact on bounded domains. Thus
generally the embedding of $\WLPhispace$ into  Orlicz space
$\LPhispace[\Phi_n]$ is not compact but with our assumptions on $N$ we get
$\WLPhispace(\Omega)\hookrightarrow\hookrightarrow\LPhispace[N](\Omega)$.

For technical reasons we define
\[\widehat{\Phi_n}: \R^n\to[0,\infty],\quad\widehat{\Phi_n}(v)=\Phi_n(|v|).\]
and assume \eqref{as:Phi<< Phinzdachem}.
In isotropic case assumption \eqref{as:Phi<< Phinzdachem} is trivial but it 
occurs to be nontrivial in anisotropic case. Even when $\Phi$ is a 
sum of powers we have situations where $\Phi$ grows faster in some direction 
than $\Phi_n$.
\begin{example}\label{ex:growth}

	\begin{enumerate}
		\item Consider $\Phi:\R^4\to[0,\infty)$,  
$\Phi(v)=|v_1|^2+|v_2|^2+|v_3|^2+|v_4|^7$. Then $\Phi_n(|v|)\approx |v|^{56/9}.$ 
Since $\frac{56}{9}<7$, assumption \eqref{as:Phi<< Phinzdachem} is not 
satisfied.
		\item Consider  $\Phi:\R^3\to[0,\infty)$,  
$\Phi(v)=|v_1|^2+|v_2|^2+|v_3|^7$. Then $\Phi_n(|v|)\approx|v|^{21}$. Since 
$21>7$, assumption \eqref{as:Phi<< Phinzdachem} is satisfied.
	\end{enumerate}
\end{example}

Function $\widehat{\Phi_n}$ and assumption \eqref{as:Phi<< Phinzdachem} will 
be important in the next lemma which will show the connection between vector 
valued functions and scalar functions. It will be useful when we transform
scalar valued function $u$ into vector valued $u\xi$ which will become useful
in applying H\"{o}lder inequality in the proof of lemma \ref{lem:PSconvergeAE}. 
We generalize a well known fact that in the isotropic Orlicz-Sobolev spaces on 
bounded domains weak convergence in $\WLPhispace$ implies convergence in 
$\LPhispace$. 

\begin{lemma}\label{lem:Phinzdachem}
	Assume that \eqref{as:Phi<< Phinzdachem} holds and $u_n\rightharpoonup u$ in 
	$\WLPhispace(\Omega)$, where $\Omega\subset\R^n$ is such that  
$|\Omega|<\infty$. Let $\xi\in\mathbf{L}^{\infty}(\Omega)$ be a vector 
	valued function. Then (up to a subsequence) $u_n\xi\to u\xi$ in 
	$\mathbf{L}^{\Phi}(\Omega)$.
\end{lemma} 
\begin{proof}
	Due to \cite{Cia00} we have 
	$\WLPhispace(\Omega)\embed\mathbf{L}^{\Phi_n}(\Omega)$. Notice that by the 
	definition we have $u_n\in \mathbf{L}^{\Phi_n}(\Omega)\Rightarrow u_n\xi\in 
	\mathbf{L}^{\widehat{\Phi_n}}(\Omega)$. By \eqref{as:Phi<< Phinzdachem} we 
know that 
	$\mathbf{L}^{\widehat{\Phi_n}}(\Omega)\hookrightarrow\mathbf{L}^{\Phi}
	(\Omega)$ compactly, so there exists a subsequence $u_n\xi\to u\xi$ in 
	$\mathbf{L}^{\Phi}(\Omega)$.
\end{proof}

We will use indices describing the growth of $\Phi$ and $N$:
\begin{definition}\label{def:indices}
\[
i_{\Phi}=\inf_{|v|>0}\frac{v\nabla \Phi(v)}{\Phi(v)}, \quad
s_{\Phi}=\sup_{|v|>0}\frac{v\nabla\Phi(v)}{\Phi(v)}\quad \text{where 
$v\in\R^n$},\]
\[
i_{N}=\inf_{t>0}\frac{tN'(t)}{N(t)}, \quad                    
s_{N}=\sup_{t>0}\frac{tN'(t)}{N(t)},\]
\end{definition}

Notice that by \eqref{N:growth} function $N$ has the same growth as 
$\Phi_{\circ}$ so it also satisfies $\Delta_2,\,\nabla_2$ conditions. As a 
consequence $1<i_{\Phi}\leq s_{\Phi}<\infty$ and $1<i_N\leq s_N<\infty$. Let
$\underline{\xi_{A}}(t)=\min\{t^{i_{A}},t^{s_{A}}\}$, 
$\overline{\xi_{A}}(t)=\max\{t^{i_{A}},t^{s_{A}}\}$ for all $t\geq 0$.

\begin{lemma}[cf. lem. 2.3 in \cite{AlvFigGioSan14}, Appendix A in 
\cite{ChmMak19}]\label{lem:PhimodularNormcon}
If $\Phi$ is n-dimensional N-function, and $N$ is N-function  which satisfy 
$\Delta_2,\,\nabla_2$ globally, then:
	\[
	\underline{\xi_{\Phi}}(\LPhinorm{u}{\R^n})\leq
	\int_{\R^n}\Phi(u)\leq \overline{\xi_{\Phi}}(\LPhinorm{u}{\R^n})
	\]
	for $u\in\LPhispace(\R^n)$.
	
		\[
	\underline{\xi_{N}}(\LPhinorm[N]{u}{\R^n})\leq
	\int_{\R^n}N(|u|)\leq \overline{\xi_{N}}(\LPhinorm[N]{u}{\R^n})
	\]
	for $u\in\LPhispace[N](\R^n)$.
\end{lemma}

Next lemma shows relation between pointwise and weak convergence in Orlicz 
spaces. It essentially is the lemma 1.4 from \cite{Gos74} rewritten in our 
notation. Notice 
that in \cite{Gos74} the lemma is formulated and proved for isotropic 
N-function but 
examination of the proof shows that it can be replaced by 
G-function (with vector valued $u_n$). Later in the proof of lemma 
\ref{lem:critical} 
we will use both versions.

\begin{lemma}\label{lem:weak_conv}
	Assume that $u_n$ is a sequence in $\mathbf{L}^{M}(\mathbb{R}^n)$ for some 
	N-function (or G-function) $M$ satisfying conditions  $\Delta_2, \nabla_2$
and $u_n\to u$ a.e. in 
	$\mathbb{R}^n$. Assume also that $M(u_n)\leq w_n$ a.e. in $\mathbb{R}^n$ and 
	$w_n\to w$ in $\mathbf{L}^1$. Then $u\in \mathbf{L}^{M}(\mathbb{R}^n)$ and 
	$u_n\rightharpoonup u$ in $\mathbf{L}^{M}(\mathbb{R}^n)$.
\end{lemma}

The following theorem comes from \cite{AlvFigGioSan14}. We formulate it 
for function $N$. Such modification can be made because by 
\eqref{N:growth} the Sobolev conjugate function 
$\Phi_n$ defined for $\Phi$ is equivalent to the conjugate function of $N$.
\begin{theorem}[Lions type lemma]\label{thm:lions}
	Let $\{u_n\}$ be a bounded sequence in the space $\WLPhispace(\R^n)$
and 
\begin{equation}\label{eq:lionslimit}
\lim_{k\to\infty}\left[\sup_{y\in\R^n}\int_{B_r(y)}N(|u_k|)\right]=0 
\quad\text{ for some $r>0$}.
\end{equation}
Assume also that $M:\R\to[0,\infty)$ is a Young 
function satisfying $\Delta_2$ globally and

	\begin{equation}\label{eq:phi/psi=0inzero}\tag{$M_1$}
	 	\lim_{t\to 0} \frac{M(t)}{N(t)}=0,
	\end{equation}
	\begin{equation}\label{eq:phi/psi=0ininfty}\tag{$M_2$}
	\lim_{t\to \infty} \frac{M(t)}{\Phi_n(t)}=0,
\end{equation}
Then \[
u_n \to 0 \quad \text{in} \quad\LPhispace[M](\R^n).
\]
	\end{theorem}

The following technical lemma is inspired by lemma 6 in \cite{DalMaMur98}. 

\begin{lemma}\label{lem:monotone}
	Assume that $\Phi$ is a differentiable, strictly convex G-function. If for 
some $\{ v_n\}\subset\WLPhispace(\R^n)$ 
$$(\nabla\Phi(\nabla v_n)-\nabla\Phi(\nabla v))(\nabla 
v_n-\nabla v)\to 0 \text{ a.e.}$$ then also  $\nabla v_n\to\nabla v$ a.e.
\end{lemma}
\begin{proof}
	On the contrary assume that there exists a subsequence $\{ v_n\}$ such that 
$|\nabla v_n-\nabla 
	v|\geq \delta>0$ on some set $\Omega$ of positive measure. Define 
	$$t_n=\frac{\delta}{|\nabla v_n-\nabla v|}, \quad \xi_n=t_n \nabla 
	v_n+(1-t_n)\nabla v.$$
	It is straightforward to check that 
	$$\xi_n-\nabla v = t_n(\nabla v_n-\nabla v), \quad \nabla v_n-\xi_n = 
	(1-t_n)(\nabla v_n-\nabla v),$$ 
	$$\quad |\xi_n-\nabla v|=\delta, \quad 0<t_n\leq 
	1.$$
	
	From the equations above we get
	\begin{align*}
		&(\nabla\Phi(\xi_n)-\nabla\Phi(\nabla v))(\nabla v_n-\nabla v)\geq 
		(\nabla\Phi(\xi_n)-\nabla\Phi(\nabla v))(\xi_n-\nabla v)>0,\\
		&(\nabla\Phi(\nabla v_n)-\nabla\Phi(\xi_n))(\nabla v_n-\nabla v)\geq 
		(\nabla\Phi(\nabla v_n)-\nabla\Phi(\xi_n))(\nabla v_n-\xi_n)> 0
	\end{align*}
	and after adding the inequalities
	\begin{align*}
		&(\nabla\Phi(\nabla v_n)-\nabla\Phi(\nabla v))(\nabla v_n-\nabla v)= \\
		&(\nabla\Phi(\xi_n)-\nabla\Phi(\nabla v))(\nabla v_n-\nabla v)+
		(\nabla\Phi(\nabla v_n)-\nabla\Phi(\xi_n))(\nabla v_n-\nabla v)\geq\\
		&(\nabla\Phi(\xi_n)-\nabla\Phi(\nabla v))(\nabla v_n-\nabla v)\geq 
		(\nabla\Phi(\xi_n)-\nabla\Phi(\nabla v))(\xi_n-\nabla v)>0.
	\end{align*}
	As we assumed $(\nabla\Phi(\nabla v_n)-\nabla\Phi(\nabla v))(\nabla 
v_n-\nabla 	v)\to 0$ so by the obtained inequality also 
	$$(\nabla\Phi(\xi_n)-\nabla\Phi(\nabla v))(\xi_n-\nabla v)\to 0 
\quad\text{a.e. in }\Omega$$ 
	For a.e. $x\in\Omega$ we have $|\xi_n(x)-\nabla v(x)|=\delta$ so sequence 
$\xi_n(x)$ is bounded in $\mathbb{R}^n$. Thus we can choose a convergent 
	subsequence: $\xi_n(x)\to\xi$ and immediately $|\xi-\nabla v(x)|=\delta$ 
(so $\xi\neq\nabla v$). Recall that $\nabla \Phi$ is continuous so 
	$$(\nabla\Phi(\xi_n)-\nabla\Phi(\nabla v))(\xi_n-\nabla v)\to 
	(\nabla\Phi(\xi)-\nabla\Phi(\nabla v))(\xi-\nabla v)= 0$$
	which is a contradiction with $\xi\neq\nabla v$.
\end{proof}

\section{Mountain Pass Geometry}

Let us denote by 
$\JSymbol:\WLPhispace(\R^n)\to\R$ the energy functional related to 
\eqref{eq:AQP} given by 
\begin{equation}\label{eq:Jot} \JSymbol(u)=\int_{\R^N}\Phi(\nabla 
u)+V(x)N(u)-F(u)\,dx. \end{equation}
It is easy to check (see \cite{BarCia17}), that 
$\JSymbol\in C^1(\WLPhispace(\R^n),\R)$  with
\begin{equation}\label{eq:J} \JSymbol'(u)v=\int_{\R^N}\nabla\Phi(\nabla 
u)\nabla v+V(x)N'(u)v-f(u)v\,dx \end{equation} 
for all $v\in\WLPhispace$ and 
that critical points of $\JSymbol$ are weak solutions of \eqref{eq:AQP}. 
The following lemmas show that $\JSymbol$ (under all the previous assumptions) 
has the mountain pass geometry.
\begin{lemma}
	There exist $\alpha, \rho>0$, such that $\JSymbol(u)\geq\alpha$, if  
$\WLPhinorm{u}{\R^n}=\rho$.
\end{lemma}
\begin{proof}
	Choose $\rho>0$ and $u$, such that $$\WLPhinorm{u}{\R^n}=\rho<1 \text{ and 
} \left[\int_{\R^n}\Phi(\nabla u)\,dx\right]^{1/n}\leq \frac{1}{K},$$where 
$K$ is 
a constant from the inequality \eqref{eq:poinc:int}.
	
	From assumptions \eqref{f:in0}-\eqref{f:takitrocheAR}, given  
$\varepsilon>0$, there exists $C_{\varepsilon}>0$ such that 
	\begin{equation}
		\theta F(t)\leq f(t)t\leq \varepsilon N(t)+C_{\varepsilon}\Phi_n(t) 
\quad \text{ for all }t\geq 0.
	\end{equation}
	Hence
	\[
	\JSymbol(u)\geq\int_{\R^N}\Phi(\nabla u)+N(u)\left(V(x)- \frac{ 
\varepsilon}{\theta}\right)-\frac{C_{\varepsilon}}{\theta}\Phi_n(|u|)\,dx
	\]
	By the inequality \eqref{eq:poinc:int} and convexity of $\Phi_n$ we obtain
	\[
	\JSymbol(u)\geq\int_{\R^N}\Phi(\nabla u)+N(u)\left(V(x)- \frac{ 
\varepsilon}{\theta}\right)\,dx-\frac{KC_{\varepsilon}}{\theta}\left(\int_{\R_n}
\Phi(\nabla u)\,dx\right)^{1+1/n}
	\]
	Of course we could have chosen $\varepsilon$ as small enough that $V(x)-
\frac{\varepsilon}{\theta}$ is positive and bounded away from $0$. Now for some
positive constants we can write
	\[
	\JSymbol(u)\geq C_1\left(\int_{\R^N}\Phi(\nabla
u)+N(u)\,dx\right)-C_2\left(\int_{\R^N}\Phi(\nabla
u)+N(u)\,dx\right)^{\frac{n+1}{n}}
	\]

	Let us define a function $h(t)=C_1 t-C_2 t^{\frac{n+1}{n}}$. By
calculating zeros of $h'$ it is easy to check that $h$ is positive and
increasing for $0<t<t_0=\left(\frac{n}{n+1}\frac{C_1}{C_2}\right)^n$. Functions
$\underline{\xi_{\Phi}}$ and $\underline{\xi_{N}}$ are increasing therefore we
can chose $\rho\leq 1$ such that
$\underline{\xi_{\Phi}}(\rho)+\underline{\xi_{N}}(\rho)\leq t_0$. Now it is
easy to see that if $\WLPhinorm{u}{\R^n}=\rho$ then by lemma
\ref{lem:PhimodularNormcon}
	\begin{align*}
	 \int_{\R^N}\Phi(\nabla u)+N(u)\,dx&\leq
\overline{\xi_{\Phi}}(\LPhinorm{\nabla
u}{\R^n})+\overline{\xi_{N}}(\LPhinorm[N]{u}{\R^n})\\
&\leq \overline{\xi_{\Phi}}(\rho)+\overline{\xi_{N}}(\rho)\leq t_0.
	\end{align*}
	Also by lemma \ref{lem:PhimodularNormcon} we have
	\[
	\underline{\xi_{\Phi}}(\LPhinorm{\nabla
u}{\R^n})+\underline{\xi_{N}}(\LPhinorm[N]{u}{\R^n})\leq \int_{\R^N}\Phi(\nabla
u)+N(u)\,dx
	\]
	Let $\beta=\LPhinorm[N]{u}{\R^n}$, then
$$\underline{\xi_{\Phi}}(\LPhinorm{\nabla
u}{\R^n})+\underline{\xi_{N}}(\LPhinorm[N]{u}{\R^n})
=\underline{\xi_{\Phi}}(\rho-\beta)+\underline{\xi_{N}}(\beta)\quad\text{and}
\quad \rho , \beta\in[0,1]$$
	therefore by the definitions of $\underline{\xi_{\Phi}}$ and
$\underline{\xi_{N}}$:
$$\underline{\xi_{\Phi}}(\rho-\beta)+\underline{\xi_{N}}(\beta)=(\rho-\beta)^{s_
{ \Phi}}+\beta^{s_N}.$$ It can be easily verified that
$\min\limits_{\beta\in[0,\rho]}((\rho-\beta)^{s_{
\Phi}}+\beta^{s_N})=\gamma_0>0$, thus also $\int_{\R^N}\Phi(\nabla
u)+N(u)\,dx\geq\gamma_0$.

By the properties of $h$ we see that $\JSymbol(u)\geq
h(\gamma_0)$.
\end{proof}

\begin{lemma}
\[
\JSymbol(tv)\to-\infty \text{ as  }t\to\infty \quad\text{ for any }  v\in 
C_0^{\infty}(\R^n)
\]
\end{lemma}

\begin{proof}
	From assumption \eqref{f:takitrocheAR}, there exists $C_5,C_6>0$ such that 
	$F(t)\geq C_5|t|^{\theta}-C_6$ for all $t\in\R\backslash\{0\}.$
	Set  $v\in C_0^{\infty}(\R^n)\backslash\{0\}$. Then, from lemma  
\ref{lem:PhimodularNormcon},  we have
\begin{align*}
	\JSymbol(tv)&=\int_{\R^N}\Phi(t\nabla v)+V(x)N(tv)-F(tv)\,dx\\
	&\leq   
\overline{\xi_{\Phi}}\left(\LPhinorm{tv}{\R^n}\right)+\overline{\xi_{N}}
(\LPhinorm[N]{tv}{\R^n})-C_5\int_{\R^n}|tv|^{\theta}dx+C_6\lambda(\mathrm
{ supp } \, v).
\end{align*}
	Hence, for sufficiently large $t$, we obtain
	\begin{align*}
		\JSymbol(tv)&\leq t^{s_{\Phi}}\overline{\xi_{\Phi}}
\left(\LPhinorm{v}{\R^n}\right)+ t^{s_N}
\overline{\xi_{N}}(\LPhinorm[N]{v}{\R^n})\\
&-C_5t^{\theta}\LPhinorm[1]{v}{\R^n}
+C_6\lambda(\mathrm{supp}\, v).
	\end{align*}
 Since $s_{\Phi}<\theta$ and $s_N<\theta$, the result follows.
\end{proof}

Now we can apply a version of the mountain pass theorem without the  
Palais-Smale condition (see \cite{BreNir91}) to get a sequence $\{u_n\}$, which 
satisfies
\begin{equation}\label{eq:PS_sequence}
\JSymbol(u_n)\to c\quad\text{ and }\JSymbol'(u_n)\to 0  \quad \text{as }  
n\to\infty,
\end{equation}
where the level $c$ is characterized by
\[
c=\inf_{\gamma\in\Gamma}\max_{s\in[0,1]}\JSymbol(\gamma(s))>0,
\]
and \[\Gamma=\{\gamma\in C([0,1],\WLPhispace(\R^n))\colon \JSymbol(0)=0\, \text{ 
and }\,\JSymbol(\gamma(1))<0\}.
\]

Note that we obtain only a minimizing sequence with a property that if there 
exists a limit then it is a critical point. The convergence of the sequence 
defined in \eqref{eq:PS_sequence} is not obvious in contrast to the classical 
mountain pass theorem. It will be examined in the next chapter.

\section{Convergence}

In this chapter we will prove several lemmas concerning the convergence of 
the sequence defined in \eqref{eq:PS_sequence}. 
\begin{lemma}\label{lem:PSbounded} If $\{v_n\}$ is the sequence 
defined in \eqref{eq:PS_sequence}, then $\{v_n\}$ is bounded in 
$\WLPhispace(\R^n)$. 
\end{lemma}

\begin{proof}
 We shall start with finding upper bounds for 
 \begin{align*}
&\JSymbol(v_n)-\frac{1}{\theta}\JSymbol'(v_n)v_n =\int_{\R^n}
\Phi(\nabla v_n)+V(x)N(v_n)-F(v_n)\,dx\\
&-\frac{1}{\theta}
\int_{\R^n}\nabla\Phi(\nabla v_n)\nabla 
v_n+V(x)N'(v_n)v_n-f(v_n)v_n\,dx \\&=\frac{1}{\theta}\int_{\R^n}
\left(\theta\Phi(\nabla v_n)-\nabla\Phi(\nabla v_n)\nabla v_n\right) 
+V(x)\left( \theta N(v_n)- N'(v_n)v_n \right) \\
&+ \left(f(v_n)v_n 
-\theta F(v_n)\right)\,dx.
\end{align*}
By \eqref{f:takitrocheAR} we have $f(v_n)v_n 
-\theta F(v_n)\geq 0$. Similarly, since $s_NN(t)\geq t N'(t)$ and 
$s_{\Phi}\Phi(\xi)\geq\xi\nabla\Phi(\xi)$ we have 
$$\theta N(v_n)- N'(v_n)v_n= 
(\theta-s_N)N(v_n)+s_N N(v_n)-N'(v_n)v_n\geq (\theta-s_N)N(v_n)$$
\begin{align*}
\theta\Phi(\nabla v_n)-\nabla\Phi(\nabla v_n)\nabla 
v_n&=(\theta-s_{\Phi})\Phi(\nabla v_n)+s_{\Phi} \Phi(\nabla v_n)- 
\nabla\Phi(\nabla v_n)\nabla v_n\\
&\geq (\theta-s_{\Phi})\Phi(\nabla v_n)
\end{align*}
Those inequalities combined with the formula for $\JSymbol'$ lead to
\begin{align*}
&\int_{\R^n} \frac{\theta-s_{\Phi}}{\theta}\Phi(\nabla
v_n)+\frac{\theta-s_{N}}{\theta}V(x)N(v_n)\,dx \\
&\leq \JSymbol(v_n)-\frac{1}{\theta}\JSymbol'(v_n)v_n\leq
c(1+\|v_n\|_{\WLPhispace}).
\end{align*}
 Now by lemma \ref{lem:PhimodularNormcon} and \eqref{as:V>0} we get 
$$\underline{\xi_{\Phi}}(\LPhinorm{\nabla v_n}{\R^n})+V_0 
\underline{\xi_{N}}(\LPhinorm[N]{v_n}{\R^n})
\leq\int_{\R^n} \Phi(\nabla v_n)+V(x)N(v_n)\,dx$$
and as a conclusion
$$\underline{\xi_{\Phi}}(\LPhinorm{\nabla v_n}{\R^n})+V_0 
\underline{\xi_{N}}(\LPhinorm[N]{v_n}{\R^n})\leq
c(1+\|\nabla v_n\|_{\LPhispace(\R^n)}+\|v_n\|_{\LPhispace[N](\R^n)}).$$
For large $\LPhinorm{\nabla v_n}{\R^n}$ the inequality can not 
hold because on the left-hand side there are powers of $\LPhinorm{\nabla 
v_n}{\R^n}$ and $\|v_n\|_{\LPhispace[N](\R^n)}$ but on the right-hand side 
there are only $\LPhinorm{\nabla v_n}{\R^n}$ and 
$\|v_n\|_{\LPhispace[N](\R^n)}$. Thus existence of a subsequence of 
$\{v_n\}$ unbounded in $\WLPhispace$ occurs to be impossible.
\end{proof}

Notice that, by reflexivity of the space and lemma \ref{lem:PSbounded}, there 
exists a subsequence of $\{v_n\}$ (still denoted by $\{v_n\}$) weakly convergent 
to some $v\in\WLPhispace(\R^n)$.

\begin{lemma}\label{lem:PSconvergeAE} $\nabla v_n \to \nabla v$ almost 
everywhere in $\R^n$. \end{lemma}

\begin{proof} By convexity of $\Phi$ we have 
\begin{equation}\label{eq:PhiMonotone}
 \left(\nabla\Phi(x)-\nabla\Phi(y)\right)(x-y)\geq 0, \quad \forall x,y\in\R^n.
\end{equation}
For some fixed $R>0$ we choose $\xi_R\in C^{\infty} (\R^n)$ such that 
$0\leq\xi_R(x)\leq 1$, $\xi_R(x)=1$ for $x\in B_R(0)$ and 
$\textrm{supp}\,\xi_R\subset 
B_{2R}(0)$. Definition of $\xi_R$ together with \eqref{eq:PhiMonotone} gives us 
\begin{align}\label{eq:ineq_conv_Phi}
 0\leq&\int\limits_{B_R(0)} (\nabla\Phi(\nabla v_n)-\nabla\Phi(\nabla 
v))(\nabla v_n-\nabla 
v) dx\nonumber\\
\leq& \int\limits_{B_{2R}(0)} (\nabla\Phi(\nabla 
v_n)-\nabla\Phi(\nabla 
v))(\nabla 
v_n-\nabla v)\xi_R dx\\
=& \int\limits_{B_{2R}(0)} \nabla\Phi(\nabla v_n)(\nabla v_n-\nabla v)\xi_R 
dx-\int\limits_{B_{2R}(0)} \nabla\Phi(\nabla v)(\nabla v_n-\nabla v)\xi_R dx 
\nonumber
\end{align}
Sequence $\{(v_n-v)\xi_R\}$ is bounded in $\WLPhispace(\R^n)$. Indeed
$\nabla((v_n-v)\xi_R)=\nabla(v_n-v)\xi_R+(v_n-v)\nabla\xi_R$, the first
summand is bounded in $\LPhispace[\Phi]$, the second is bounded in
$\LPhispace[N]$ thus (thanks to \eqref{as:Phi<< Phinzdachem}) also bounded in
$\LPhispace[\Phi]$.

Since $\JSymbol'(v_n)\to0$:
\begin{align}\label{eq:Jprim_to_0}
 &\int\limits_{B_{2R}(0)} \nabla\Phi(\nabla v_n)\nabla((v_n-v)\xi_R)dx+ 
 \int\limits_{B_{2R}(0)} V(x)N'(|v_n|)(v_n-v)\xi_R dx\nonumber\\
 &-\int\limits_{B_{2R}(0)} f(v_n)(v_n-v)\xi_R dx \to 0.
\end{align}
Sequence $\{N'(|v_n|)\}$ is bounded in 
$\mathbf{L}^{\widetilde{N}}(B_{2R}(0))$ by lemma  \ref{lem:PSbounded} and
$$\int\limits_{B_{2R}(0)} \widetilde{N}(N'(v_n))dx\leq 
\int\limits_{B_{2R}(0)} N(2v_n) dx \leq 
K \int\limits_{B_{2R}(0)} N(v_n)dx.$$
Together with the boundedness of $V$ (which is a result of its periodicity) and 
H\"{o}lder inequality we get
\begin{align*}
 &\left| \int\limits_{B_{2R}(0)} V(x)N'(|v_n|)(v_n-v)\xi_R dx \right|\leq
 \int\limits_{B_{2R}(0)} |V(x)|N'(|v_n|)|v_n-v| dx\\
 &\leq 2M \| N'(|v_n|) \|_{\mathbf{L}^{\widetilde{N}}(B_{2R}(0))}\| v_n-v 
\|_{\mathbf{L}^{N}(B_{2R}(0))}
\leq c_1 \| v_n-v \|_{\mathbf{L}^{N}(B_{2R}(0))} \to 0
\end{align*}
since $\WLPhispace\hookrightarrow\hookrightarrow \LPhispace[N]$. From 
continuity of $f$, assumption \eqref{f:ininfty},  the H\"older inequality and 
properties of the conjugate of the Young function (see (2.9) in 
\cite{AlChlCiaZat19}) for arbitrarily small $\epsilon$ we obtain

\begin{align*}
	&\left|\int\limits_{B_{2R}(0)} f(v_n)(v_n-v)\xi_R dx\right| \leq 
	\int\limits_{B_{2R}(0)} \left(c_2+\epsilon 
	\frac{\Phi_n(|v_n|)}{v_n}\right)|v_n-v|dx \\
	&\leq  c_3 \| v_n-v \|_{\mathbf{L}^{N}(B_{2R}(0))}+2\epsilon \|	
\Phi_n(|v_n|)/v_n
	\|_{\mathbf{L}^{\Phi_n^{\ast}}(B_{2R}(0))}\| v_n-v
	\|_{\mathbf{L}^{{\Phi}_n}(B_{2R}(0))}\\
	&\leq  c_3 \| v_n-v \|_{\mathbf{L}^{N}(B_{2R}(0))}+2\epsilon \|v_n
	\|_{\mathbf{L}^{\Phi_n}(B_{2R}(0))}\| v_n-v
	\|_{\mathbf{L}^{{\Phi}_n}(B_{2R}(0))}.
\end{align*}

By theorem \ref{thm:embcian} the sequence $\{ v_n\}$ is bounded 
in $\LPhispace[\Phi_n](B_{2R}(0))$ and convergent $\mathbf{L}^{N}(B_{2R}(0))$ 
as $N\prec\prec \Phi_n$. Thus
$$\left|\int\limits_{B_{2R}(0)} f(v_n)(v_n-v)\xi_R dx\right|\to0.$$

Of course
\begin{align}\label{eq:grad_of_product}
 &\int\limits_{B_{2R}(0)} \nabla\Phi(\nabla 
v_n)\nabla((v_n-v)\xi_R)dx\nonumber\\
 &= \int\limits_{B_{2R}(0)} \nabla\Phi(\nabla v_n)\nabla(v_n-v)\xi_Rdx+ 
 \int\limits_{B_{2R}(0)} \nabla\Phi(\nabla v_n)(v_n-v)\nabla\xi_Rdx
\end{align}
In the next step we consider the last integral
$$\int\limits_{B_{2R}(0)} \nabla\Phi(\nabla v_n)(v_n-v)\nabla\xi_R\, dx.$$
We shall use similar inequalities as before:
\begin{align*}
&\left|\int\limits_{B_{2R}(0)} \nabla\Phi(\nabla 
v_n)(v_n-v)\nabla\xi_R\,dx\right|\\
&\leq 2\| \nabla\Phi(\nabla v_n) 
\|_{\mathbf{L}^{\widetilde{\Phi}}(B_{2R}(0))} 
\| (v_n-v)\nabla\xi_R\|_{\mathbf{L}^{\Phi}(B_{2R}(0))}.
\end{align*}
Of course
\begin{equation}\label{eq:bounded_Phi_tilde}
\int\limits_{B_{2R}(0)} \widetilde{\Phi}(\nabla\Phi(v_n))dx\leq
\int\limits_{B_{2R}(0)} \Phi(2v_n) dx \leq 
K \int\limits_{B_{2R}(0)} \Phi(v_n)dx<\infty
\end{equation}
thus $\| \nabla\Phi(\nabla v_n) \|_{\mathbf{L}^{\widetilde{\Phi}}(B_{2R}(0))}
$ is bounded.

The convergence $\|
(v_n-v)\nabla\xi_R\|_{\mathbf{L}^{\Phi}(B_{2R}(0))}\to
0$ is a result of lemma \ref{lem:Phinzdachem} and again we have
$$\int\limits_{B_{2R}(0)} \nabla\Phi(\nabla v_n)(v_n-v)\nabla\xi_R\, dx\to 0.$$

Now from \eqref{eq:Jprim_to_0} and \eqref{eq:grad_of_product} we get 
\begin{equation}\label{eq:conv_nabla_v_n}
\int\limits_{B_{2R}(0)} \nabla\Phi(\nabla v_n)(\nabla v_n-\nabla v)\xi_R\, 
dx\to 0.
\end{equation}
It is easy to see that $\nabla\Phi(\nabla v)\xi_R\in 
\LPhispace[\widetilde{\Phi}](B_{2R}(0))$. From weak convergence 
$v_n\rightharpoonup v$ in $\WLPhispace(B_{2R}(0))$
$$\int\limits_{B_{2R}(0)} \nabla\Phi(\nabla v)\xi_R\nabla v_n\, dx\to 
\int\limits_{B_{2R}(0)} \nabla\Phi(\nabla v)\xi_R\nabla v\,dx.$$

Applying the above and \eqref{eq:conv_nabla_v_n} to \eqref{eq:ineq_conv_Phi} 
leads us to
$$\int\limits_{B_R(0)} (\nabla\Phi(\nabla v_n)-\nabla\Phi(\nabla 
v))(\nabla v_n-\nabla v) dx\to 0$$
so we can choose a subsequence still denoted by $v_n$ such that 
$$(\nabla\Phi(\nabla v_n)-\nabla\Phi(\nabla 
v))(\nabla v_n-\nabla v)\to 0, \,\text{a.e. in}\, B_R(0).$$

By lemma \ref{lem:monotone} $\nabla v_n\to\nabla v$ a.e. in $B_R(0)$.
\end{proof}

Now we will prove that $v$ is a critical point for $\JSymbol$. 
\begin{lemma}\label{lem:critical}
 Let $X=\mathbf{W}^1_0\mathbf{L}^{\Phi}(\mathbb{R}^n)$ be the subspace of 
functions with compact support. For all $u\in X$ we have $\JSymbol'(v)u=0$.
\end{lemma}
\begin{proof}
 In the proof of lemma \ref{lem:PSconvergeAE} in \eqref{eq:bounded_Phi_tilde}
we have shown that the sequence $\nabla\Phi(\nabla v_n)$ is bounded in
$\LPhispace[\widetilde{\Phi}](B_{2R}(0))$ but of course the same inequality
could be written with $\R^n$ instead of a ball. By lemma \ref{lem:weak_conv}
$$\int\limits_{\mathbf{R}^n}\nabla\Phi(\nabla v_n)\nabla u \, dx\to
\int\limits_{\mathbf{R}^n}\nabla\Phi(\nabla v)\nabla u \, dx.$$
Analogously by the boundedness of $V(x)N'(|v_n|)$ in
$\mathbf{L}^{\widetilde{N}}$ and $f(v_n)$ in $\mathbf{L}^{\widetilde{\Phi_n}}$ 
we see that 
\begin{align*}
\int\limits_{\mathbf{R}^n}V(x)N'(|v_n|)u\, dx&\to 
\int\limits_{\mathbf{R}^n}V(x)N'(|v|)u\, dx \quad\text{and}\\ 
\int\limits_{\mathbf{R}^n}f(|v_n|)u\, dx&\to 
\int\limits_{\mathbf{R}^n}f(|v|)u\, dx.
\end{align*}
\end{proof}

Finally $v$ is a critical point of $\JSymbol$ by the density of $X$ in 
$\mathbf{W}^1\mathbf{L}^{\Phi}(\mathbb{R}^n)$. To prove the main theorem 
\ref{thm:main} we need to check that $v\neq 0$. To do that we will use the 
Lions type lemma (theorem \ref{thm:lions}). 

\begin{lemma}
 There exists $R>0$ such that 
 $$\liminf\limits_{n\to\infty}\sup\limits_{y\in\mathbb{R}^n}\, \int_{B_R(y)} 
N(v_n)\, dx>0.$$
\end{lemma}
\begin{proof}
 If there would not exist any such $R$ then by theorem \ref{thm:lions} the 
following 
$$\int_{\mathbb{R}^n} M(v_n)\, dx\to 0$$
would hold for any function $M$ satisfying \eqref{eq:phi/psi=0inzero} and 
\eqref{eq:phi/psi=0ininfty}.  By assumptions \eqref{f:in0} and \eqref{f:ininfty}
for some fixed $\epsilon>0$ we have $|f(v_n)v_n|\leq \epsilon N(|v_n|)$ and 
$|f(v_n)v_n|\leq\epsilon \Phi_n(|v_n|)$ in neighbourhoods of $0$ and $\infty$. 
Outside those neighbourhoods for some function $M$ satisfying 
\eqref{eq:phi/psi=0inzero} and \eqref{eq:phi/psi=0ininfty} there exists 
$C_{\epsilon}$ such that $|f(v_n)v_n|\leq C_{\epsilon}M(v_n)$ so on  whole $\R$ 
we have
$$|f(v_n)v_n|\leq \epsilon N(|v_n|)+ \epsilon \Phi_n(|v_n|)+ 
C_{\epsilon}M(v_n).$$
Now, as the $\epsilon$ is arbitrarily small, it is easy to see that 
$$\int_{\mathbb{R}^n} f(v_n)v_n\, dx \to 0.$$
Notice that $\JSymbol'(v_n)v_n\to 0$, so
$$\int_{\R^N}\nabla\Phi(\nabla v_n)\nabla v_n+V(x)N'(v_n)v_n\, dx \to 0.$$
Since $\Phi, N\in\Delta_2,\nabla_2$ we have $i_{\Phi}\Phi(\nabla v_n)\leq 
\nabla\Phi(v_n)\nabla v_n$ and analogous inequality for $N$. Thus
$$\int_{\R^N}\Phi(\nabla v_n)+V(x)N(v_n)\, dx \to 0,$$
which means that $v_n\to 0$ in $\WLPhispace(\mathbb{R}^n)$. That would imply 
$c=0$ in \eqref{eq:PS_sequence} which is impossible.
\end{proof}

\section{Finish of the proof of the main theorem}
As a result of the previous 
lemma there exist $R,q>0$ and $\{ y_n \}\subset 
\mathbb{Z}^n$ such that $\int_{B_R(y_n)}N(v_n)\, dx > q$. Define 
$w_n(x)=v_n(x-y_n)$. By periodicity of $V$: $\JSymbol(w_n)=\JSymbol(v_n)$, 
$\JSymbol'(w_n)\to0$ and $\|w_n\|=\|v_n\|$. As before, there exists $w$ such 
that $w_n\rightharpoonup w$. The embedding 
$\WLPhispace(B_R(0))
\hookrightarrow\mathbf{L}^{N}(B_R(0))$ is compact, so   
$\int_{B_R(y_n)}N(w_n)\, dx > q>0$ implies $\int_{B_R(y_n)}N(w)\, dx\geq 
q>0$. Thus $w\neq 0$.

\vskip 50pt
\noindent\textbf{Funding.} There is no funding to declare.\\
\noindent\textbf{Author contribution.} Author declares that he has written the 
entire manuscript.\\
\textbf{Conflict of interest.} The author declares that he has no conflict 
of interests.

\bibliographystyle{abbrv}
\bibliography{biblio}


\end{document}